\documentclass[12pt]{article}
\usepackage{amsfonts,amscd,a4}
\usepackage{color}
\usepackage{graphicx}
\usepackage{theorem}
\newtheorem{theo}{Theorem}[section]
\newtheorem{prop}[theo]{Proposition}
\newtheorem{lemma}[theo]{Lemma}

{\theorembodyfont{\rm}

}
\newcommand{\bA}{{\bf A}}

\newcommand{\bI}{{\bf I}}

\newcommand{\cC}{{\mathcal C}}

\newcommand{\cF}{{\mathcal F}}
\newcommand{\cG}{{\mathcal G}}

\newcommand{\cJ}{{\mathcal J}}

\newcommand{\cS}{{\mathcal S}}

\newcommand{\cY}{{\mathcal Y}}

\newcommand{\eA}{{\sf A}}
\newcommand{\eB}{{\sf B}}

\newcommand{\eD}{{\sf D}}

\newcommand{\eO}{{\sf O}}

\newcommand{\eT}{{\sf T}}

\newcommand{\sC}{{\mathbb C}}

\newcommand{\sF}{{\mathbb F}}

\newcommand{\sN}{{\mathbb N}}

\newcommand{\sZ}{{\mathbb Z}}
\newcommand{\qed}{\rule{1ex}{1ex}}

\newcommand{\im}{\mbox{\rm im} \,}

\newcommand{\op}{\mbox{\rm Op} \,}

\begin{document}
\title{Finite sections of band-dominated operators on discrete groups}
\author{Vladimir S. Rabinovich, Steffen Roch}
\date{}
\maketitle
\begin{abstract}
Let $\Gamma$ be a finitely generated discrete exact group. We consider operators on $l^2(\Gamma)$ which are composed by operators of multiplication by a function in $l^\infty (\Gamma)$ and by the operators of left-shift by elements of $\Gamma$. These operators generate a $C^*$-subalgebra of $L(l^2(\Gamma))$ the elements of which we call band-dominated operators on $\Gamma$. We study the stability of the finite sections method for band-dominated operators with respect to a given generating system of $\Gamma$. Our approach is based on the equivalence of the stability of a sequence and the Fredholmness of an associated operator, and on Roe's criterion for the Fredholmness of a band-dominated operator on a exact discrete group, which we formulate in terms of limit operators. Special emphasis is paid to the quasicommutator ideal of the algebra generated by the finite sections sequences and to the stability of sequences in that algebra. For both problems, the sequence of the discrete boundaries plays an essential role.
\end{abstract}
\section{Introduction} \label{s1}
Let $\Gamma$ be a countable (not necessarily commutative) discrete group. We write the group operation as multiplication and let $e$ stand for the identity element of $\Gamma$. For each non-empty subset $X$ of $\Gamma$, let $l^2(X)$ stand for the Hilbert space of all functions $f : X \to \sC$ with
\[
\|f\|^2 := \sum_{x \in X} |f(x)|^2 < \infty.
\]
For $X = \emptyset$, we define $l^2(X)$ as the space $\{0\}$ consisting of the zero element only. We consider $l^2(X)$ as a closed subspace of $l^2(\Gamma)$ in a natural way. The orthogonal projection from $l^2(\Gamma)$ to $l^2(X)$ will be denoted by $P_X$. Thus, $P_\Gamma$ and $P_\emptyset$ are the identity and the zero operator, respectively. For $s \in \Gamma$, let $\delta_s$ be the function on $\Gamma$ which is 1 at $s$ and 0 at all other points. The family $(\delta_s)_{s \in \Gamma}$ forms an orthonormal basis of $l^2(\Gamma)$, to which we refer as the standard basis.

The {\em left regular representation} $L : \Gamma \to L(l^2(\Gamma))$ of $\Gamma$ associates with every group element $r$ a unitary operator $L_r$ such that $L_r \delta_s = \delta_{rs}$ for $s \in \Gamma$. Since $\delta_{rs}(t) = \delta_s(r^{-1}t)$, one has $(L_r u)(t) = u(r^{-1}t)$ for every $u \in l^2(\Gamma)$. Hence, $r \mapsto L_r$ is a group isomorphism. Further, we associate with each function $a \in l^\infty(\Gamma)$ the operator $aI$ of multiplication by $a$, i.e., $(au)(t) = a(t) u(t)$ for $u \in l^2(\Gamma)$. The smallest closed subalgebra of $L(l^2(\Gamma))$ which contains all operators $L_r$ with $r \in \Gamma$ and $aI$ with $a \in l^\infty(\Gamma)$ is called the {\em algebra of the band-dominated operators on} $\Gamma$. We denote it by $\eB\eD\eO(\Gamma)$. Besides $\eB\eD\eO(\Gamma)$ we consider the smallest closed subalgebra ${\sf Sh}(\Gamma)$ of $L(l^2(\Gamma))$ which contains all "shift" operators $L_r$ with $r \in \Gamma$. Clearly, the algebras $\eB\eD\eO(\Gamma)$ and ${\sf Sh}(\Gamma)$ are symmetric and, hence, $C^*$-subalgebras of $L(l^2(\Gamma))$.

Let $\cY = (Y_n)_{n = 1}^\infty$ be an increasing sequence of finite subsets of $\Gamma$ with $\cup_{n \ge 1} Y_n = \Gamma$. A sequence $(A_n)_{n = 1}^\infty$ of operators $A_n : \im P_{Y_n} \to \im P_{Y_n}$ is called {\em stable} if there is an $n_0 \ge 1$ such that the operators $A_n$ are invertible for $n \ge n_0$ and the norms of their inverses $A_n^{-1}$ are bounded uniformly with respect to $n \ge n_0$. Note that stability is crucial for many questions in asymptotic numerical analysis. It dominates topics like the approximate solution of operator equations and the approximate spectral and pseudo-spectral theory. For a detailed overview see \cite{HRS2}.

Let $A \in L(l^2(\Gamma))$. The operators $P_{Y_n} A P_{Y_n} : \im P_{Y_n} \to \im P_{Y_n}$ are called the {\em finite sections} of $A$ with respect to $\cY$. In this paper, we are interested in the stability of the finite sections sequence $(P_{Y_n} A P_{Y_n})$ when $A \in \eB\eD\eO(\Gamma)$. The finite sections method for band-dominated operators on the group $\sZ$ of the integers is quite well understood, see \cite{RRS2,RRSB,RRS6,Roc10}. Finite sections for operators in ${\sf Sh}(\Gamma)$ with an arbitrary exact countable discrete group $\Gamma$ were considered in \cite{Roc12}.

Our approach to study the stability of the finite sections method for operators in $\eB\eD\eO(\Gamma)$ is close to that in \cite{Roc10,Roc12}. We make use of the fact that a sequence $(A_n)$ is stable if and only if an associated operator has the Fredholm property. In case the $A_n$ are the finite sections of a band-dominated operator, the associated operator is a band-dominated operator again. So the desired stability result will finally follow from Roe's criterion for the Fredholm property of band-dominated operators in \cite{Roe2}. We thus start with recalling Roe's result in Section \ref{s2}.

In Section \ref{c61}, we provide an algebraic frame to study the stability of operator sequences. We introduce the $C^*$-algebra $\cS_\cY ({\sf BDO} (\Gamma))$ generated by all finite sections sequences $(P_{Y_n} A P_{Y_n})$ with $A \in \eB\eD\eO(\Gamma)$ and show that this algebra splits into the direct sum of ${\sf BDO} (\Gamma)$ and of an ideal which can be characterized as the quasicommutator ideal of the algebra. A main result is that the sequence $(P_{\partial Y_n})$ of the discrete boundaries always belongs to the algebra $\cS_\cY ({\sf BDO} (\Gamma))$, and that this sequence already generates the quasicommutator ideal. This surprising fact has been already observed in other settings, for example for the algebras $\cS(\eT(C))$ of the finite sections method for the Toeplitz operators (a classical result, closely related to the present paper) and $\cS_\cY ({\sf Sh} (\Gamma))$ (see \cite{Lin1} for the group $\Gamma = \sZ^n$ and \cite{Roc12} for the general case), but also for the finite sections algebra $\cS(\eO_N)$ related with a concrete representation $\eO_N$ of the Cuntz algebra (see \cite{Roc11}).

The final Section 4 is devoted to the prove of the stability theorem. We employ Roe's criterion using the limit operators language from \cite{RRSB}. The main task is to compute all (or at least a sufficient number of) limit operators of the band-dominated operator associated with a finite sections sequence.

The work on this paper was supported by CONACYT Project 81615 and DFG Grant Ro 1100/8-1.
\section{The algebra of the band-dominated operators} \label{s2}
We start with some alternate characterizations of band-dominated operators and the algebra generated by them. Consider functions $k \in l^\infty (\Gamma \times \Gamma)$ with the property that there is a finite subset $\Gamma_0$ of $\Gamma$ such that $k(t, \, s) = 0$ whenever $t s^{-1} \not\in \Gamma_0$. Then
\begin{equation} \label{1.1groups}
(Au)(t) := \sum_{s \in \Gamma} k (t, \, s) \, u(s), \qquad t \in \Gamma,
\end{equation}
defines a linear operator $A$ on the linear space of all functions $u : \Gamma \to \sC$, since the occurring series is finite for every $t \in G$. We call operators of this form {\em band operators} and the set $\Gamma_0$ a {\em band-width} of $A$.
\begin{prop}
A operator in $L(l^2(\Gamma))$ is a band operator if and only if it can be written as a finite sum $\sum b_i L_{t_i}$ where $b_i \in l^\infty(\Gamma)$ and $t_i \in \Gamma$.
\end{prop}
{\bf Proof.} Let $A$ be an operator of the form (\ref{1.1groups}) and let $\Gamma_0 := \{t_1, \, t_2, \, \ldots, \, t_r \}$ be a finite subset of $\Gamma$ such that $k(t, \, s) = 0$ if $t s^{-1} \not\in \Gamma_0$ (or, equivalently, if $s$ is not of the form $t_i^{-1} t$ for some $i$). Thus,
\[
(Au)(t) = \sum_{i = 1}^r k (t, \, t_i^{-1} t) \, u(t_i^{-1} t) \quad \mbox{for} \; t \in \Gamma.
\]
Set $b_i(t) := k (t, \, t_i^{-1} t)$. The functions $b_i$ are in $l^\infty (\Gamma)$, and one has
\begin{equation} \label{bdo}
A = \sum_{i = 1}^r b_i L_{t_i}.
\end{equation}
Conversely, one easily checks that each operator $L_t$ with $t \in \Gamma$ is a band operator with band width $\{t\}$ and that each operator $bI$ with $b \in l^\infty(\Gamma)$ is a band operator with band width $\{e\}$. Since the band operators form an algebra, each finite sum $\sum b_i L_{t_i}$ is a band operator. \hfill \qed \\[3mm]
It is easy to see that the representation of a band operator on $\Gamma$ in the form (\ref{bdo}) with $b_i \neq 0$ is unique. The functions $b_i$ are called the {\em diagonals} of the operator $A$. In particular, operators in ${\sf Sh}(\Gamma)$ can be considered as band-dominated operators with constant coefficients.

It is easy to see that the band operators form a symmetric algebra of bounded operators on $l^2(\Gamma)$. The norm closure of that algebra is just the algebra ${\sf BDO}(\Gamma)$, and this is why we call the elements of that algebra {\em band-dominated} operators.

The algebras ${\sf BDO}(\Gamma)$ and ${\sf Sh}(\Gamma)$ occur at many places and under different names in the literature. The algebra ${\sf Sh}(\Gamma)$ is $^*$-isomorphic to the {\em reduced group $C^*$-algebra} $C^*_r(\Gamma)$ in a natural way (see Section 2.5 in \cite{BroO}). It can thus be considered as a concrete faithful representation of $C^*_r(\Gamma)$. Note also that the reduced group $C^*$-algebra coincides with the universal group $C^*$-algebra $C^*(\Gamma)$ if the group $\Gamma$ is amenable. For this and further characterizations of amenable groups, see Theorem 2.6.8 in \cite{BroO}. The algebra ${\sf BDO}(\Gamma)$ occurs in coarse geometry and is known there as the {\em uniform Roe algebra} or the {\em reduced translation algebra} (\cite{Roe1}). It can be identified with the reduced crossed product of the $C^*$-algebra $l^\infty(\Gamma)$ with the group $\Gamma$ when the group action $\alpha : \Gamma \to \mbox{Aut} \, l^\infty(\Gamma)$ is specified as
\[
(\alpha_g f)(t) := f(g^{-1} t)
\]
for $f \in l^\infty(\Gamma)$ and $g, \, t \in \Gamma$. Note that amenability of $\Gamma$ is not needed for the following result. But if $\Gamma$ is amenable, then the reduced crossed product $l^\infty(\Gamma) \times_{\alpha r} \Gamma$ coincides with the full crossed product $l^\infty(\Gamma) \times_\alpha \Gamma$ (see \cite{Ped1}, Theorem 7.7.7 and \cite{Dav1}, Corollary VII.2.2).
\begin{theo}
The reduced crossed product $l^\infty(\Gamma) \times_{\alpha r} \Gamma$ of the $C^*$-dynamical system $(l^\infty(\Gamma), \, \Gamma, \, \alpha)$ is $^*$-isomorphic to ${\sf BDO}(\Gamma)$.
\end{theo}
{\bf Proof.} Let $l^2(\Gamma, l^2(\Gamma))$ stand for the Hilbert space of all functions $x : \Gamma \to l^2(\Gamma)$ with $\sum_{s \in \Gamma} \|x(s)\|^2 < \infty$. For $a \in l^\infty(\Gamma)$, let $\pi(a)$ denote the operator $aI$ of multiplication by $a$ on $l^2 (\Gamma)$ and define an operator $\tilde{\pi} (a)$ on $l^2(\Gamma, l^2(\Gamma))$ by
\[
(\tilde{\pi} (a) x)(s) := \pi ( \alpha_s^{-1} (a)) (x(s)).
\]
For $g \in \Gamma$, let $\tilde{L}_g$ be the operator on $l^2(\Gamma, l^2(\Gamma))$ defined by
\[
(\tilde{L}_g x)(s) := x(t^{-1} s).
\]
The pair $(\tilde{\pi}, \, \tilde{L})$ constitutes a covariant representation of the $C^*$-dynamical system $(l^\infty(\Gamma), \, \Gamma, \, \alpha)$ on $l^2(\Gamma, l^2(\Gamma))$. By the definition of the reduced crossed product (see \cite{Bla1,Dav1,Ped1}, for instance), $l^\infty(\Gamma) \times_{\alpha r} \Gamma$ is the smallest $C^*$-subalgebra of $L(l^2(\Gamma, l^2(\Gamma)))$ which contains all operators $\tilde{\pi} (a)$ and $\tilde{L}_g$ with $a \in l^\infty(\Gamma)$ and $g \in \Gamma$. One can show (\cite{Ped1}, Theorem 7.7.5) that each faithful representation $(\pi^\prime, \, H)$ of $l^\infty(\Gamma)$ in place of the representation $(\pi, \, l^2(\Gamma))$ leads to the same algebra.

We identify $l^2(\Gamma, \, l^2(\Gamma))$ with $l^2(\Gamma \times \Gamma)$ via the mappings
\[
J : l^2(\Gamma, \, l^2(\Gamma)) \to l^2(\Gamma \times \Gamma), \; (Jx)(s, \,n) := (x(s))(n),
\]
\[
J^{-1} : l^2(\Gamma \times \Gamma) \to l^2(\Gamma, \, l^2(\Gamma)), \; ((J^{-1}y)(s))(n) := y(s, \, n)
\]
and determine the corresponding operators
\[
\hat{\pi} (a) := J \tilde{\pi}(a) J^{-1} \quad \mbox{and} \quad \hat{L}_g := J \tilde{L}_g J^{-1}.
\]
A straightforward calculation gives
\begin{equation} \label{ea1.7}
(\hat{\pi} (a) x)(s, \, n) = a(sn) x(s, \, n) \quad \mbox{and} \quad (\hat{L}_g x) (s, \, n) = x(g^{-1} s, \, n).
\end{equation}
Let $\cC$ refer to the smallest $C^*$-subalgebra of $L(l^2(\Gamma \times \Gamma))$ which contains all operators
$\hat{\pi} (a)$ and $\hat{L}_g$ with $a \in l^\infty(\Gamma)$ and $g \in \Gamma$, given by (\ref{ea1.7}). For $n \in \Gamma$, let
\[
H_n := \{ x \in l^2(\Gamma \times \Gamma): x(s, \, m) = 0 \; \mbox{whenever} \; m \neq n \}.
\]
We identify $l^2(\Gamma \times \Gamma)$ with the orthogonal sum $\oplus_{n \in \Gamma} H_n$ such that $x \in l^2(\Gamma \times \Gamma)$ is identified with $\oplus h_n \in \oplus H_n$ where $h_n(s) = x(s, \, n)$. From (\ref{ea1.7}) we conclude that each space $H_n$ is invariant with respect to each operator in $\cC$ (i.e., $A H_n \subseteq H_n$ for $A \in \cC$). Hence, each operator $A \in \cC$ corresponds to a diagonal matrix operator $\mbox{diag} \, ( \ldots, \, A_n, \, A_{n+1}, \, \ldots)$ with respect to the decomposition of $l^2(\Gamma \times \Gamma)$ into the orthogonal sum of its subspaces $H_n$. Thus, $A_n$ is the restriction of $A$ onto $H_n$.

Let $\cC_n$ be the $C^*$-algebra of all restrictions of operators in $\cC$ onto $H_n$. It is clear that each of the spaces $H_n$ is isometric to $l^2(\Gamma)$, with the isometry given by
\[
J_n : H_n \to l^2(\Gamma), \quad (J_n x)(s) := x(s, \, n),
\]
\[
J_n^{-1} : l^2(\Gamma)  \to H_n, \quad (J_n^{-1} x)(s, \, n) := x(s).
\]
Then
\begin{eqnarray*}
(J_n \hat{\pi}(a) J_n^{-1} x)(s) & = & (\hat{\pi}(a) J_n^{-1} x)(s, \, n) \; = \; (a(sn) (J^{-1}x)) (s, \, n) \\
& = & a(sn) x(s) \; = \; (R_n \pi(a) R_n^{-1} x)(s)
\end{eqnarray*}
where $(R_n f)(s) = f(sn)$ stands for the operator of the right-regular representation of $\Gamma$. Similarly,
\begin{eqnarray*}
(J_n \hat{L}_g J_n^{-1} x)(s) & = & (\hat{L}_g J_n^{-1} x)(s, \, n) \; = \; (J^{-1} x)(g^{-1}s, \, n) \\
& = & x(g^{-1} s) \; = \; (L_g x)(s).
\end{eqnarray*}
Thus,
\[
J_n \hat{\pi}(a) J_n^{-1} = R_n \pi(a) R_n^{-1} \quad \mbox{and} \quad  J_n \hat{L}_g J_n^{-1} = L_g = R_n L_g R_n^{-1}.
\]
Consequently, the mapping
\[
{\sf BDO}(\Gamma) \to \cC, \quad A \mapsto \mbox{diag} \, ( \ldots, \, J_n^{-1} R_n A R_n^{-1} J_n, \, \ldots)
\]
is a $^*$-isomorphism. Since $\cC$ is evidently $^*$-isomorphic to the reduced crossed product $l^\infty(\Gamma) \times_{\alpha r} \Gamma$, the assertion follows. \hfill \qed \\[3mm]
Our next goal is to recall Roe's criterion \cite{Roe2} for the Fredholm property of band-dominated operators on $l^2(\Gamma)$. We are going to formulate this criterion in the language of limit operators.

Let $h : \sN \to \Gamma$ be a sequence tending to infinity in the sense that for each finite subset $\Gamma_0$ of $\Gamma$, there is an $n_0 \in \sN$ such that $h(n) \not\in \Gamma_0$ if $n \ge n_0$. Clearly, if $h$ tends to infinity, then the inverse sequence $h^{-1}$ tends to infinity, too. We say that an operator $A_h \in L(l^2(\Gamma))$ is a {\em limit operator of $A \in L(l^2(\Gamma))$ defined by the sequence} $h$ if
\[
R_{h(m)}^{-1} A R_{h(m)} \to A_h \quad \mbox{and} \quad R_{h(m)}^{-1} A^* R_{h(m)} \to A_h^*
\]
strongly as $m \to \infty$ (as before, the $R_r$ are given by the right-regular representation of $\Gamma$ on $l^2(\Gamma)$). Clearly, every operator has at most one limit operator with respect to a given sequence $h$. Note that the generating function of the shifted operator $R_r^{-1} A R_r$ is related with the generating function of $A$ by
\begin{equation} \label{equa0}
k_{R_r^{-1} A R_r} (t, \, s) = k_A (t r^{-1}, \, s r^{-1})
\end{equation}
and that the generating functions of $R_{h(m)}^{-1} A R_{h(m)}$ converge pointwise on $\Gamma \times \Gamma$ to the generating function of the limit operator $A_h$ (if the latter exists).

It is an important property of band-dominated operators that they always possess limit operators. More general, the following result can be proved by a standard Cantor diagonal argument (see \cite{RRS1,RRS2,RRSB}).
\begin{prop} \label{prop1}
Let $A$ be a band-dominated operator on $l^2(\Gamma)$. Then every sequence $h : \sN \to \Gamma$ which tends to infinity possesses a subsequence $g$ such that the limit operator $A_g$ of $A$ with respect to $g$ exists.
\end{prop}
Let $A$ be a band-dominated operator and $h : \sN \to \Gamma$ a sequence tending to infinity for which the limit operator $A_h$ of $A$ exists. Let $B$ be another band-dominated operator. By Proposition \ref{prop1} we can choose a subsequence $g$ of $h$ such that the limit operator $B_g$ exists. Then the limit operators of $A$, $A+B$ and $AB$ with respect to $g$ exist, and
\[
A_g = A_h, \qquad (A+B)_g = A_g + B_g, \qquad (AB)_g = A_g B_g.
\]
Thus, the mapping $A \mapsto A_h$ acts, at least partially, as an algebra homomorphism.

The following theorem is due to Roe \cite{Roe2}, see also \cite{RaR5}. Recall that a group $\Gamma$ is called {\em exact}, if its reduced translation algebra is exact as a $C^*$-algebra. The latter algebra is defined as the reduced crossed product of $l^\infty(\Gamma)$ by $\Gamma$ and coincides with the $C^*$-algebra of all band-dominated operators on $l^2(\Gamma)$ in our setting. The class of exact groups is extremely rich. It includes all amenable groups (hence, all solvable groups such as the discrete Heisenberg group and the commutative groups) and all hyperbolic groups (in particular, all free groups with finitely many generators) (see \cite{Roe1}, Chapter 3).
\begin{theo}[Roe] \label{t1.1}
Let $\Gamma$ be a finitely generated discrete and exact group, and let $A$ be a band-dominated operator on $l^2(\Gamma)$. Then the operator $A$ is Fredholm on $l^2(\Gamma)$ if and only if all limit operators of $A$ are invertible and if the norms of their inverses are uniformly bounded.
\end{theo}
Note that this result holds as well if the left regular representation is replaced by the right regular one and if, thus, the operators $L_s$ and $R_t$ change their roles. In fact, the results of \cite{RaR5,Roe2} are presented in this symmetric setting. In \cite{RaR5} we showed moreover that the uniform boundedness condition in Theorem \ref{t1.1} is redundant for band operators if the group $\Gamma$ has sub-exponential growth and if not every element of $\Gamma$ is cyclic in the sense that $w^n = e$ for some positive integer $n$. For details see \cite{RaR5}. Note that the condition of sub-exponential growth is satisfied by the abelian groups $\sZ^N$, the discrete Heisenberg group and, more general, by nilpotent groups (in fact, these groups have polynomial growth), whereas the growth of the free group $\sF_N$ with $N > 1$ is exponential.
\begin{theo} \label{t120309.1}
Let $\Gamma$ be a finitely generated discrete and exact group with sub-exponential growth which possesses at least one non-cyclic element, and let $A$ be a band operator on $l^2(\Gamma)$. Then the operator $A$ is Fredholm on $l^2(\Gamma)$ if and only if all limit operators of $A$ are invertible.
\end{theo}
\section{The algebra of the finite sections method} \label{c61}
Given an increasing sequence $\cY := (Y_n)_{n \ge 1}$ of finite subsets of $\Gamma$ such that $\cup_{n \ge 1} Y_n = \Gamma$, let $\cF_\cY$ denote the set of all bounded sequences $\bA = (A_n)$ of operators $A_n : \im P_{Y_n} \to \im P_{Y_n}$. Equipped with the operations
\[
(A_n) + (B_n) := (A_n + B_n), \quad
(A_n) (B_n) := (A_n B_n), \quad (A_n)^* := (A_n^*)
\]
and the norm
\[
\|\bA\|_{\cF_\cY} := \|A_n\|,
\]
the set $\cF_\cY$ becomes a $C^*$-algebra with identity $\bI = (Y_n)$, and the set $\cG_\cY$ of all sequences $(A_n) \in \cF_\cY$ with $\lim \|A_n\| = 0$ forms a closed ideal of $\cF_\cY$. The relevance of the algebra $\cF_\cY$ and its ideal $\cG_\cY$ in our context stems from the fact (following by a simple Neumann series argument) that a sequence $\bA \in \cF_\cY$ is stable if, and only if, its coset $\bA + \cG_\cY$ is invertible in the quotient algebra $\cF_\cY/\cG_\cY$. Thus, every stability problem is equivalent to an invertibility problem in a suitably chosen $C^*$-algebra.

Let further stand $\cF^C_\cY$ for the set of all sequences $\bA = (A_n)$ of operators $A_n : \im P_{Y_n} \to \im P_{Y_n}$ with the property that the sequences $(A_n P_{Y_n})$ and $(A_n^* P_{Y_n})$ converge strongly. By the uniform boundedness principle, the quantity $\sup \|A_n P_{Y_n}\|$ is finite for every sequence $(A_n)$ in $\cF^C_\cY$. Thus, $\cF^C_\cY$ is a closed and symmetric subalgebra of $\cF_\cY$ which contains $\cG_\cY$, and the mapping
\begin{equation} \label{e91.5}
W : \cF^C_\cY \to L(l^2(X)), \quad (A_n) \mapsto \mbox{s-lim} \, A_n P_{Y_n}
\end{equation}
is a $^*$-homomorphism. Note that $\bI \in \cF^C_\cY$ and that $W(\bI)$ is the identity operator $I$ on $L^2(\Gamma)$.

For each $C^*$-subalgebra $\eA$ of $L(l^2(\Gamma))$, write $D$ for the mapping of finite sections (or spatial) discretization, i.e.,
\begin{equation} \label{e91.6}
D : L(l^2(\gamma)) \to \cF_\cY, \quad A \mapsto (P_{Y_n} A P_{Y_n}),
\end{equation}
and let $\cS_\cY (\eA)$ stand for the smallest closed $C^*$-subalgebra of the algebra $\cF_\cY$ which contains all sequences $D(A)$ with $A \in \eA$. Clearly, $\cS_\cY (\eA)$ is contained in $\cF^C_\cY$, and the mapping $W$ in (\ref{e91.5}) induces a $^*$-homomorphism from $\cS_\cY (\eA)$ onto $\eA$. On this level, one cannot say much about the algebra $\cS_\cY (\eA)$. The simple proof of the following is in \cite{Roc11}.
\begin{prop} \label{p91.10}
Let $\eA$ be a $C^*$-subalgebra of $L(l^2(\Gamma))$. Then the finite sections discretization $D : \eA \to \cF_\cY$ is an isometry, and $D(\eA)$ is a closed subspace of the algebra $\cS_\cY (\eA)$. This algebra splits into the direct sum
\[
\cS_\cY(\eA) = D(\eA) \oplus (\ker W \cap \cS_\cY(\eA)),
\]
and for every operator $A \in \eA$ one has
\[
\|D(A)\| = \min_{K \in \ker W} \|D(A) + K\|.
\]
Finally, $\ker W \cap \cS_\cY(\eA)$ is equal to the quasicommutator ideal of $\cS_\cY (\eA)$, i.e., to the smallest closed ideal of $\cS_\cY (\eA)$ which contains all sequences
$(P_{Y_n} A_1 P_{Y_n} A_2 P_{Y_n} - P_{Y_n} A_1 A_2 P_{Y_n})$ with operators $A_1, \, A_2 \in \eA$.
\end{prop}
We denote the ideal $\ker W \cap \cS_\cY(\eA)$ by $\cJ_\cY (\eA)$. Since the first item in the decomposition $D(\eA) \oplus \cJ_\cY(\eA)$ of $\cS_\cY(\eA)$ is isomorphic (as a linear space) to $\eA$, a main part of the description of the algebra $\cS_\cY(\eA)$ is to identify the ideal $\cJ_\cY(\eA)$.

We are going to present two alternate descriptions of the quasicommutator ideal $\cJ_\cY({\sf BDO}(\Gamma))$ of the finite sections algebra $\cS_\cY({\sf BDO}(\Gamma))$. For we have to introduce some notions of topological type. Note that the standard topology on $\Gamma$ is the discrete one; so every subset of $\Gamma$ is open with respect to this topology.

Let $\Omega$ be a finite subset of $\Gamma$ which contains the identity element $e$ and which generates $\Gamma$ as a semi-group, i.e., if we set $\Omega_0 := \{e\}$ and if we let $\Omega_n$ denote the set of all words of length at most $n$ with letters in $\Omega$ for $n \ge 1$, then $\cup_{n \ge 0} \Omega_n = \Gamma$. Note also that the sequence $(\Omega_n)$ is increasing; so the operators $P_{\Omega_n}$ can play the role of the finite sections projections $P_{Y_n}$, and in fact we will obtain some of the subsequent results exactly for this sequence.

With respect to $\Omega$, we define the following "algebro-topological" notions. Let $A \subseteq \Gamma$. A point $a  \in A$ is called an {\em $\Omega$-inner} point of $A$ if $\Omega a := \{ \omega a : \omega \in \Omega \} \subseteq A$. The set $\mbox{int}_\Omega A$ of all $\Omega$-inner points of $A$ is called the {\em $\Omega$-interior} of $A$, and the set $\partial_\Omega A := A \setminus \mbox{int}_\Omega A$ is the {\em $\Omega$-boundary} of $A$. Note that we consider the $\Omega$-boundary of a set always as a part of that set. (In this point, the present definition of a boundary differs from other definitions in the literature; see \cite{Ada1} for instance.)

One easily checks that the $\Omega$-interior and the $\Omega$-boundary of a set are invariant with respect to multiplication from the right-hand side:
\[
(\mbox{int}_\Omega A) s = \mbox{int}_\Omega (As) \quad \mbox{and} \quad (\partial_\Omega A) s = \partial_\Omega (As)
\]
for $s \in \Gamma$. One also has
\begin{equation} \label{e030309.5}
\Omega_{n-1} \subseteq \mbox{int}_\Omega \Omega_n \subseteq \Omega_n \quad \mbox{for each} \; n \ge 1,
\end{equation}
whence
\begin{equation} \label{e030309.6}
\partial_\Omega \Omega_n \subseteq \Omega_n \setminus \Omega_{n-1} \quad \mbox{for each} \; n \ge 1.
\end{equation}
Here is a first result which describes $\cJ_\cY({\sf BDO}(\Gamma))$ in terms of generators of $\Gamma$. Abbreviate $I - P_A =: Q_A$.
\begin{theo} \label{t030309.8}
$\cJ_\cY({\sf BDO}(\Gamma))$ is the smallest closed ideal of $\cS_\cY ({\sf BDO}(\Gamma))$ which contains all sequences
\begin{equation} \label{e030309.9}
(P_{Y_n} L_{\omega^{-1}} Q_{Y_n} L_\omega P_{Y_n})_{n \ge 1} \quad \mbox{with} \quad \omega \in \Omega.
\end{equation}
\end{theo}
We call $(P_{\partial_\Omega Y_n})_{n \ge 1}$ the {\em sequence of the discrete boundaries} of the finite section method with respect to $(Y_n)$. Note that the assumptions in the following theorem are satisfied if $Y_n = \Omega_n$ due to (\ref{e030309.5}).
\begin{theo} \label{t030309.14}
Assume that $Y_{n-1} \subseteq \mbox{\rm int}_\Omega Y_n \subseteq Y_n$ for all $n \ge 2$ and that $\cup_{n \ge 1} Y_n = \Gamma$. Then the sequence $(P_{\partial_\Omega Y_n})_{n \ge 1}$ of the discrete boundaries belongs to the algebra $\cS_\cY ({\sf BDO}(\Gamma))$, and the quasicommutator ideal is generated by this sequence, i.e., $\cJ_\cY ({\sf BDO}(\Gamma))$ is the smallest closed ideal of $\cS_\cY ({\sf BDO}(\Gamma))$ which contains $(P_{\partial_\Omega Y_n})_{n \ge 1}$.
\end{theo}
Both results were proved in \cite{Roc12} for the ideal $\cJ_\cY ({\sf Sh}(\Gamma))$ of $\cS_\cY ({\sf Sh}(\Gamma))$ in place of $\cJ_\cY ({\sf BDO}(\Gamma))$. The above theorems follow from these results since every multiplication operator $aI$ commutes with every projection $P_Y$ where $Y \subseteq \Gamma$.
\section{Stability} \label{c63}
We are now going to study the stability of sequences in $\cS_\cY ({\sf BDO}(\Gamma))$ via the limit operators method.
The key observations are that the stability of a sequence in that algebra is equivalent to the Fredholm property of a certain associated operator, which is band-dominated, such that the Fredholm property of that operator can be studied by means of its limit operators via Roe's result.

Let again $\cY := (Y_n)$ be an increasing sequence of finite subsets of $\Gamma$ with $\cup_{n \ge 1} Y_n = \Omega$. A sequence $(v_n) \subseteq \Gamma$ is called an {\em inflating} sequence for $\cY$ if $Y_m v_m^{-1} \cap Y_n v_n^{-1} = \emptyset$ for $m \neq n$. The existence of inflating sequences is easy to see. Moreover, the following lemma was shown in \cite{Roc12}.
\begin{lemma}
Let $\cY = (Y_n)$ be as above and $V$ an infinite subset of $\Gamma$. Then there is an inflating sequence for $\cY$ in $V$.
\end{lemma}
In what follows we choose and fix an inflating sequence $(v_n)$ for $\cY$ and set 
\begin{equation} \label{e180209.2}
\Gamma^\prime := \Gamma \setminus \cup_{n=1}^\infty Y_n v_n^{-1}.
\end{equation}
For $s \in \Gamma$, let again $R_s : l^2(\Gamma) \to l^2(\Gamma)$ refer to the operator $(R_s f)(t) := f(ts)$. Evidently, $R_s L_t = L_t R_s$ for $s, \, t \in \Gamma$. The proof of the following theorem is in \cite{Roc12}.
\begin{theo} \label{t160209.6}
Let $\bA = (A_n) \in \cF_\cY$. Then \\[1mm]
$(a)$ the series
\begin{equation} \label{e160209.0}
\sum_{n=1}^\infty R_{v_n} A_n R_{v_n}^{-1}
\end{equation}
converges strongly on $l^2(\Gamma)$. The sum of this series is denoted by $\op (\bA)$. \\[1mm]
$(b)$ the sequence $(A_n)$ is stable if and only if the operator $\op (\bA) + P_{\Gamma^\prime}$ is Fredholm on $l^2(\Gamma)$. \\[1mm]
$(c)$ The mapping $\op$ is a continuous homomorphism from $\cF_\cY$ to $L(l^2(\Gamma))$.
\end{theo}
The applicability of Roe's result to the study the stability of the finite section method for band-dominated operators rests of the following fact.
\begin{prop} \label{p170209.1}
Let $\bA$ be a sequence in $\cS_\cY ({\sf BDO}(\Gamma))$. Then $\op (\bA)$ is a band-dominated operator.
\end{prop}
{\bf Proof.} First let $A \in {\sf BDO}(\Gamma)$ be a band operator and let $\Gamma_0$ be a band width of $A$. It is easy to check that then $R_{v_n} P_{Y_n} A P_{Y_n} R_{v_n}^{-1}$ is a band operator with the same band width for every $n$. The inflating property ensures that $\op \left( (P_{Y_n} A P_{Y_n}) \right)$ is a band operator with band width $\Gamma_0$, too. Now Theorem \ref{t160209.6} $(c)$ yields the assertion. \hfill \qed \\[3mm]
In order to verify the stability of a sequence $\bA \in \cS_\cY ({\sf BDO}(\Gamma))$ via the above results, we thus have to compute the limit operators of $\op (\bA) + P_{\Gamma^\prime}$, which will be our next goal. Note that the exactness of $\Gamma$ is not relevant in this computation.

Let $\Omega$ be a finite subset of $\Gamma$ with $e \in \Omega$ which generates $\Gamma$ as a semi-group and define $\Omega_n$ as above. By Theorem \ref{t160209.6}, the Fredholm property of an operator $\op(\bA)$ is independent of the concrete choice of the inflating sequence. For technical reasons, we choose an inflating sequence $(v_n)$ for the sequence
\[
\left( (Y_n \cup \Omega_n) (Y_n \cup \Omega_n)^{-1} (Y_n \cup \Omega_n) \right)_{n \ge 1}
\]
instead of $(Y_n)_{n \ge 1}$. Since
\[
Y_n \cup \Omega_n \subset (Y_n \cup \Omega_n) (Y_n \cup \Omega_n)^{-1} \subset (Y_n \cup \Omega_n) (Y_n \cup \Omega_n)^{-1} (Y_n \cup \Omega_n),
\]
$(v_n)$ is also an inflating sequence for $(Y_n)$. Moreover, since $\mbox{s-lim} \, P_{\Omega_n} = P_\Gamma = I$, one also has
\begin{equation} \label{e180209.7}
\mbox{s-lim} \, P_{(Y_n \cup \Omega_n) (Y_n \cup \Omega_n)^{-1}} = P_\Gamma = I.
\end{equation}
Let now $\bA = (A_n) \in \cS_\cY ({\sf BDO}(\Gamma))$, set as before
\[
\op(\bA) = \sum_{n=1}^\infty R_{v_n} A_n R_{v_n}^{-1} \quad \mbox{and} \quad \Gamma^\prime = \Gamma \setminus \cup_{n=1}^\infty Y_n v_n^{-1},
\]
and let $h : \sN \to \Gamma$ be a sequence tending infinity for which the limit operator
\[
(\op(\bA) + P_{\Gamma^\prime})_h := \mbox{s-lim}_{n \to \infty} R_{h(n)}^{-1} (\op(\bA) + P_{\Gamma^\prime}) R_{h(n)}
\]
exists. Then the limit operator $(\op(\bA) + P_{\Gamma^\prime})_g$ exists for every subsequence $g$ of $h$, and it coincides with $(\op(\bA) + P_{\Gamma^\prime})_h$. So we can freely pass to subsequences of $h$ if necessary. By a first passage to a suitable subsequence of $h$ we can arrange that one of the following two situations happens; so we can restrict the computation of the limit operator to these cases: \\[2mm]
\hspace*{5mm} {\bf Case 1:} All elements $h(n)$ belong to $\cup_{k \ge 1} \, v_k Y_k^{-1}$. \\[1mm]
\hspace*{5mm} {\bf Case 2:} No element $h(n)$ belongs to $\cup_{k \ge 1} \, v_k Y_k^{-1}$. \\[2mm]
We start with {\bf Case 1}. Passing again to a subsequence of $h$, if necessary, we can further suppose that each $h(n)$ belongs to one of the sets $v_k Y_k^{-1}$, say to $v_{k_n} Y_{k_n}^{-1}$, and that $v_{k_n} Y_{k_n}^{-1}$ contains no other element of the sequence $h$ besides $h(n)$. For each $n$, let $r_n$ denote the smallest non-negative integer such that $h(n) \in v_{k_n} (\partial_\Omega Y_{k_n})^{-1} \Omega_{r_n}$. Thus, $r_n$ measures the distance of $h(n)$ to the $\Omega$-boundary of $v_{k_n} Y_{k_n}^{-1}$. Set $r^* := \liminf_{n \to \infty} r_n$. Again we can distinguish two cases. \\[2mm]
{\bf Case 1.1: $r^*$ is finite.} Then there are infinitely many $n \in \sN$ such that $r_n = r^*$. Thus, there is a subsequence of $h$ (denoted by $h$ again) such that
\[
h(n) \in v_{k_n} Y_{k_n}^{-1} \cap v_{k_n} (\partial_\Omega Y_{k_n})^{-1} \Omega_{r^*} \quad \mbox{for all} \; n.
\]
Further, for each $n$ there is an $w_n^* \in \Omega_{r^*}$ such that $h(n) \in v_{k_n} (\partial_\Omega Y_{k_n})^{-1} w_n^*$. Since $\Omega_{r^*}$ is a finite set, one of the elements $w_n^*$ of $\Omega_{r^*}$ occurs for infinitely many $n$. Let $w_*$ be an element of $\Omega_{r^*}$ with this property. Consider the subsequence of $h$ which contains all elements $h(n)$ with $w_n^* = w_*$. We denote this subsequence by $h$ again and can hence assume that
\begin{equation} \label{e180209.1}
h(n) \in v_{k_n} Y_{k_n}^{-1} \cap v_{k_n} (\partial_\Omega Y_{k_n})^{-1} w_* \quad \mbox{for} \; n \ge 1.
\end{equation}
With respect to a sequence $h$ as in (\ref{e180209.1}) we obtain
\begin{eqnarray} \label{e180209.3}
\lefteqn{R_{h(n)}^{-1} (\op(\bA) + P_{\Gamma^\prime}) R_{h(n)}} \nonumber \\
&& \hspace*{-5mm} = \sum_{k=1}^\infty R_{h(n)}^{-1} R_{v_k} A_k R_{v_k}^{-1} R_{h(n)} + R_{h(n)}^{-1} P_{\Gamma^\prime} R_{h(n)} \nonumber \\
&& \hspace*{-5mm} = \sum_{k \neq k_n} R_{h(n)}^{-1} R_{v_k} A_k R_{v_k}^{-1} R_{h(n)} + R_{h(n)}^{-1} P_{\Gamma^\prime} R_{h(n)} + R_{h(n)}^{-1} R_{v_{k_n}} A_{k_n} R_{v_{k_n}}^{-1} R_{h(n)}
\end{eqnarray}
with $\Gamma^\prime$ as in (\ref{e180209.2}). By (\ref{e180209.1}), $h(n) = v_{k_n} \eta_{k_n} w_*$ with $\eta_{k_n} \in (\partial_\Omega Y_{k_n})^{-1}$. Thus, the last item in (\ref{e180209.3}) becomes
\begin{equation} \label{e180209.4}
R_{w_*^{-1}} R_{\eta_{k_n}^{-1}} A_{k_n} R_{\eta_{k_n}} R_{w_*}.
\end{equation}
Set $\Pi_n := P_{(Y_{k_n} \cup \Omega_{k_n}) (Y_{k_n} \cup \Omega_{k_n})^{-1} w_*}$. By (\ref{e180209.7}), $\Pi_n \to I$ strongly.
Since $A_{k_n}$ acts on $\im P_{Y_{k_n}}$, the operator (\ref{e180209.4}) acts on $\im P_{Y_{k_n} \eta_{k_n} w_*}$. The evident inclusion
\[
Y_{k_n} \eta_{k_n} w_* \subseteq (Y_{k_n} \cup \Omega_{k_n}) (Y_{k_n} \cup \Omega_{k_n})^{-1} w_*
\]
implies that
\[
\Pi_n R_{h(n)}^{-1} R_{v_{k_n}} A_k R_{v_{k_n}}^{-1} R_{h(n)} = R_{h(n)}^{-1} R_{v_{k_n}} A_k R_{v_{k_n}}^{-1} R_{h(n)} \Pi_n =
R_{h(n)}^{-1} R_{v_{k_n}} A_k R_{v_{k_n}}^{-1} R_{h(n)}.
\]
Let now $k \neq k_n$. Then, by the inflating property,
\begin{eqnarray} \label{e180209.8}
\lefteqn{(Y_k \cup \Omega_k) (Y_k \cup \Omega_k)^{-1} (Y_k \cup \Omega_k) v_k^{-1}} \nonumber \\
&& \cap (Y_{k_n} \cup \Omega_{k_n}) (Y_{k_n} \cup \Omega_{k_n})^{-1} (Y_{k_n} \cup \Omega_{k_n}) v_{k_n}^{-1} = \emptyset.
\end{eqnarray}
Since $Y_k v_k^{-1} \subseteq (Y_k \cup \Omega_k) (Y_k \cup \Omega_k)^{-1} (Y_k \cup \Omega_k) v_k^{-1}$ and
\[
(Y_{k_n} \cup \Omega_{k_n}) (Y_{k_n} \cup \Omega_{k_n})^{-1} \eta_{k_n}^{-1} v_{k_n}^{-1} \subseteq
(Y_{k_n} \cup \Omega_{k_n}) (Y_{k_n} \cup \Omega_{k_n})^{-1} (Y_{k_n} \cup \Omega_{k_n}) v_{k_n}^{-1}
\]
we conclude from (\ref{e180209.8}) that
\[
Y_k v_k^{-1} \cap (Y_{k_n} \cup \Omega_{k_n}) (Y_{k_n} \cup \Omega_{k_n})^{-1} \eta_{k_n}^{-1} v_{k_n}^{-1} = \emptyset
\]
whence
\[
Y_k v_k^{-1} v_{k_n} \eta_{k_n} w_* \cap (Y_{k_n} \cup \Omega_{k_n}) (Y_{k_n} \cup \Omega_{k_n})^{-1} w_* = \emptyset.
\]
Since $R_{h(n)}^{-1} R_{v_k} A_k R_{v_k}^{-1} R_{h(n)}$ is an operator living on $\im P_{Y_k v_k^{-1} v_{k_n} \eta_{k_n} w_*}$, we further conclude that
\[
R_{h(n)}^{-1} R_{v_k} A_k R_{v_k}^{-1} R_{h(n)} \Pi_n = \Pi_n R_{h(n)}^{-1} R_{v_k} A_k R_{v_k}^{-1} R_{h(n)} = 0
\]
for $k \neq k_n$. Hence,
\begin{eqnarray} \label{e180209.9}
\lefteqn{R_{h(n)}^{-1} (\op(\bA) + P_{\Gamma^\prime}) R_{h(n)}} \nonumber \\
&& = \sum_{k \neq k_n} R_{h(n)}^{-1} R_{v_k} A_k R_{v_k}^{-1} R_{h(n)} (I - \Pi_n) + R_{h(n)}^{-1} P_{\Gamma^\prime} R_{h(n)}
\nonumber \\
&& \qquad + \; R_{w_*}^{-1} R_{\eta_{k_n}}^{-1} A_{k_n} R_{\eta_{k_n}} R_{w_*} \Pi_n.
\end{eqnarray}
Since $\Pi_n \to I$ strongly, the first summand on the right-hand side of (\ref{e180209.9}) converges strongly (and even $^*$-strongly since $\Pi_n$ commutes with that sum) to zero. Thus,
\begin{eqnarray*}
\lefteqn{\mbox{s-lim} \, R_{h(n)}^{-1} (\op(\bA) + P_{\Gamma^\prime}) R_{h(n)}} \\
&& = \mbox{s-lim} \, R_{w_*}^{-1} R_{\eta_{k_n}}^{-1} A_{k_n} R_{\eta_{k_n}} R_{w_*} \Pi_n + \mbox{s-lim} \, R_{h(n)}^{-1} P_{\Gamma^\prime} R_{h(n)},
\end{eqnarray*}
provided that the strong limits on the right-hand side exist. The existence of the second strong limit can always be forced by passing to a suitable subsequence of $h$. Collecting these facts, we arrive at the following.
\begin{theo} \label{t180209.10}
Let $\bA \in \cS_\cY ({\sf BDO}(\Gamma))$, and let $h$ be a sequence such that the limit operator $\op(\bA) + P_{\Gamma^\prime}$ exists. In {\rm Case 1.1}, there is a subsequence $g$ of $h$ such that the limit operator $(P_{\Gamma^\prime})_g$ exists, and there are a monotonically increasing sequence $(k_n)$ in $\sN$, a vector $\eta_{k_n} \in (\partial_\Omega Y_{k_n})^{-1}$ for each $n \ge 1$, and a $w_* \in \Gamma$ such that
\[
(\op(\bA) + P_{\Gamma^\prime})_h = \mbox{\rm s-lim} \, R_{w_*}^{-1} R_{\eta_{k_n}}^{-1} A_{k_n} R_{\eta_{k_n}} R_{w_*} + (P_{\Gamma^\prime})_g.
\]
\end{theo}
Thus, the operator $A_{k_n}$ living on $\im P_{Y_{k_n}}$ is shifted by a vector $\eta_{k_n} \in (\partial_\Omega Y_{k_n})^{-1}$ and by another vector $w_*$ independent of $n$. It is only a matter of taste to consider $A_{k_n}$ as shifted by the vector $\eta_{k_n}^{-1}$ belonging to the $\Omega$-boundary of $Y_{k_n}$. In particular, every limit operator of $\op(\bA)$ is a shift by some vector $w_*$ of a strong limit of operators $A_{k_n}$, shifted by vectors in the $\Omega$-boundary of $Y_{k_n}$. This fact is well known for the group $\sZ$ and intervals $Y_k = [-k, \, k] \cap \sZ$ (and has been employed in \cite{RRS6} to get rid of the uniform boundedness condition in this case), and it was observed by Lindner \cite{Lin1} in case $\Gamma = \sZ^N$ and $Y_k = \Omega_k$ is a convex polygon with integer vertices.

Before turning to the other cases, let us specify Theorem \ref{t180209.10} to pure finite sections sequences for operators in ${\sf BDO}(\Gamma)$. The existence of the limit operator $(P_{\Gamma^\prime})_h$ is guaranteed if the strong limit
\begin{equation} \label{e190209.1}
\mbox{s-lim} \, R_{w_*}^{-1} R_{\eta_{k_n}}^{-1} P_{Y_{k_n}} R_{\eta_{k_n}} R_{w_*} = \mbox{s-lim} \, P_{Y_{k_n} \eta_{k_n} w_*}
\end{equation}
exists. In this case, there is a subset $\cY^{(h)}$ of $\Gamma$ such that
\begin{equation} \label{e190209.16}
\mbox{s-lim} \, P_{Y_{k_n} \eta_{k_n} w_*} = P_{\cY^{(h)}}
\end{equation}
and, thus, $(P_{\Gamma^\prime})_g = I - P_{\cY^{(h)}}$. We claim that the sequence $(\eta_{k_n} w_*)_{n \ge 1}$ tends to infinity. For this goal, it is sufficient to show that every sequence $(\mu_n)$ with $\mu_n \in \partial_\Omega Y_{k_n}$ tends to infinity. Let $\Gamma_0$ be a finite subset of $\Gamma$. Choose $n_0$ such that $\Gamma_0 \subseteq \Omega_{n_0-1}$ and $n^*$ such that $\Omega_{n_0} \subseteq Y_{k_n}$ for all $n \ge n^*$. Then $\mbox{int}_\Omega \Omega_{n_0} \subseteq \mbox{int}_\Omega Y_{k_n}$, and from (\ref{e030309.5}) we conclude that
\[
\Gamma_0 \subseteq \Omega_{n_0-1} \subseteq \mbox{int}_\Omega \Omega_{n_0} \subseteq \mbox{int}_\Omega Y_{k_n}.
\]
Hence, $\partial_\Omega Y_{k_n} \cap \Gamma_0 = \emptyset$ for all $n \ge n^*$, whence the claimed convergence.

Given a sequence $h$ such that the limit (\ref{e190209.1}) exists and a band-dominated operator $A$, let $\sigma_{op, \, h}(A)$ denote the set of all limit operators of $A$ with respect to subsequences of the sequence $(\eta_{k_n} w_*)_{n \ge 1}$. This set is not empty by Proposition \ref{prop1}.
\begin{prop}
Let $A \in {\sf BDO}(\Gamma)$, and let $h$ be a sequence such that the limit operator $\op(\bA)_h$ for the sequence $(P_{Y_n} A P_{Y_n})$ exists. In {\rm Case 1.1}, there are $k_n$, $\eta_{k_n}$ and $w_*$ as in Theorem $\ref{t180209.10}$ such that the limit $(\ref{e190209.1})$ exists. Then there is a limit operator $A_g \in \sigma_{op, \, h}(A)$ of $A$ such that
\begin{equation} \label{e190209.2}
(\op(\bA) + P_{\Gamma^\prime})_h = P_{\cY^{(h)}} A_g P_{\cY^{(h)}} + (I - P_{\cY^{(h)}}).
\end{equation}
Conversely, if the limit $(\ref{e190209.1})$ exists for a certain choice of $k_n$, $\eta_{k_n}$ and $w_*$ as in Theorem $\ref{t180209.10}$ and if $A_g$ is a limit operator of $A$ with respect to a certain subsequence $g = (\eta_{k_{n_r}} w_*)_{r \ge 1}$ of the sequence $(\eta_{k_n} w_*)_{n \ge 1}$, then the limit operator $\op(\bA)_h$ exists for the sequence $h = (v_{k_{n_r}} g_r)_{r \ge 1}$, and $(\ref{e190209.2})$ holds.
\end{prop}
{\bf Proof.} The proof of the first assertion follows easily from Theorem \ref{t180209.10}. Indeed,
\begin{eqnarray*}
\lefteqn{R_{w_*^{-1} \eta_{k_n}^{-1}} P_{Y_{k_n}} A P_{Y_{k_n}} R_{\eta_{k_n} w_*}} \\
&& = \; (R_{w_*^{-1} \eta_{k_n}^{-1}} P_{Y_{k_n}} R_{\eta_{k_n} w_*}) \cdot (R_{w_*^{-1} \eta_{k_n}^{-1}} A R_{\eta_{k_n} w_*}) \cdot (R_{w_*^{-1} \eta_{k_n}^{-1}} P_{Y_{k_n}} R_{\eta_{k_n} w_*}).
\end{eqnarray*}
The sequences in the outer parentheses converge strongly to $P_{\cY^{(h)}}$. If now $g$ is a subsequence of $(\eta_{k_n} w_*)_{n \ge 1}$ such that the limit operator $A_g$ exists, then we conclude that
\[
R_{w_*^{-1} \eta_{k_n}^{-1}} P_{Y_{k_n}} A P_{Y_{k_n}} R_{\eta_{k_n} w_*} \to P_{\cY^{(h)}} A_g P_{\cY^{(h)}}
\]
$^*$-strongly as $n \to \infty$. The second assertion is evident. \hfill \qed  \\[2mm]
{\bf Case 1.2: $r^*$ is infinite.} Recall that
\begin{equation} \label{e190209.3}
h(n) \in v_{k_n} Y_{k_n}^{-1} \quad \mbox{and} \quad
h(n) \not\in v_{k_n} (\partial_\Omega Y_{k_n})^{-1} \Omega_{r_n -1}
\end{equation}
for all $n \in \sN$. The second assertion in (\ref{e190209.3}) implies that
\[
h(n) \Omega_{r_n -1}^{-1} \cap v_{k_n} (\partial_\Omega Y_{k_n})^{-1} = \emptyset.
\]
Hence, we can rewrite (\ref{e190209.3}) as
\begin{equation} \label{e190209.4}
e \in Y_{k_n} v_{k_n}^{-1} h(n) \quad \mbox{and} \quad
\Omega_{r_n -1} \cap (\partial_\Omega Y_{k_n}) v_{k_n}^{-1} h(n) = \emptyset.
\end{equation}
We claim that this implies that
\begin{equation} \label{e190209.5}
\Omega_{r_n -1} \subseteq Y_{k_n} v_{k_n}^{-1} h(n).
\end{equation}
Suppose (\ref{e190209.5}) is wrong. Then $\Omega_{r_n -1}$ has at least one point outside $Y_{k_n} v_{k_n}^{-1} h(n)$, say $a$, but it also has points inside this set, for example the point $e$ due to the first assumption of (\ref{e190209.4}). Write $a$ as a product $a = w_{r_n-1} \ldots w_1 w_0$ of elements $w_i \in \Omega$ with $w_0 := e$, and let $0 \le j < r_n-1$ be the smallest integer such that
\[
w_j \ldots w_1 w_0 \in Y_{k_n} v_{k_n}^{-1} h(n), \quad \mbox{but} \quad  w_{j+1} w_j \ldots w_1 w_0 \not\in Y_{k_n} v_{k_n}^{-1} h(n).
\]
Then $\Omega w_j \ldots w_1 w_0 \not\subseteq Y_{k_n} v_{k_n}^{-1} h(n)$, hence
\[
w_j \ldots w_1 w_0 \in \partial_\Omega (Y_{k_n} v_{k_n}^{-1} h(n)).
\]
Since $w_j \ldots w_1 w_0 \in \Omega_{r_n -1}$, this contradicts
the second assertion of (\ref{e190209.4}), and the claim
(\ref{e190209.5}) follows. Roughly speaking, we used the fact that $\Omega$-boundaries do not have gaps. Since $P_{\Omega_n} \to I$ strongly, we conclude from (\ref{e190209.5}) that
\begin{equation} \label{e190209.8}
P_{Y_{k_n} v_{k_n}^{-1} h(n)} \to I \quad \mbox{strongly}.
\end{equation}
\begin{theo} \label{t190209.6}
Let $\bA \in \cS_\cY ({\sf BDO}(\Gamma))$ and $A := \mbox{\rm s-lim} A_n P_{Y_n}$, and let $h$ be a sequence such that the limit operator $\op(\bA)_h$ exists. Then, in {\rm Case 1.2}, either $\op(\bA)_h = R_{v^*}^{-1} A R_{v^*}$ with a fixed $v^* \in \Gamma$, or there is a limit operator $A_g$ of $A$ such that $\op(\bA)_h = A_g$. Conversely, each operator $R_{v^*}^{-1} A R_{v^*}$ with $v^* \in \Gamma$ and each limit operator $A_g$ of $A$ occur as limit operators of $\op(\bA)$.
\end{theo}
{\bf Proof.} It is sufficient to verify the assertion for pure finite sections sequences $\bA = (P_{Y_n} A P_{Y_n})$ with $A \in {\sf BDO}(\Gamma)$. For these sequences, one has
\begin{eqnarray*}
\lefteqn{R_{h(n)}^{-1} (\op(\bA) + P_{\Gamma^\prime}) R_{h(n)}} \\
&& = \; \sum_{k \neq k_n} R_{h(n)}^{-1} R_{v_k} P_{Y_k} A P_{Y_k} R_{v_k}^{-1} R_{h(n)} (I - P_{Y_{k_n} v_{k_n}^{-1} h(n)}) \\
&& \qquad + \; R_{h(n)}^{-1} P_{\Gamma^\prime} R_{h(n)} (I - P_{Y_{k_n} v_{k_n}^{-1} h(n)}) \\
&& \qquad + \; P_{Y_{k_n} v_{k_n}^{-1} h(n)} (R_{h(n)}^{-1} R_{v_{k_n}} A R_{v_{k_n}}^{-1} R_{h(n)})P_{Y_{k_n} v_{k_n}^{-1} h(n)}.
\end{eqnarray*}
Consider the sequence $(v_{k_n}^{-1} h(n))$, which is either finite or contains a subsequence which tends to infinity. In the first case, there is a $v^* \in \Gamma$ which is met by this sequence infinitely often, whence $\op(\bA)_h = R_{v^*} A R_{v^*}^{-1}$ due to (\ref{e190209.8}). In the second case, Proposition \ref{prop1} implies the existence of a subsequence $g$ of  $(v_{k_n}^{-1} h(n))$ which tends to infinity and for which the limit operator $A_g$ exists. In this case, $\op(\bA)_h = A_g$.

Conversely, given $v^* \in \Gamma$ and a limit operator $A_g$ of $A$, one can choose $h(n) := v_{k_n} v^*$ and $h(n) := v_{k_n} g(n)$ in order to obtain the limit operators $R_{v^*}^{-1} A R_{v^*}$ and $A_g$ of $\op(\bA)$, respectively. \hfill \qed \\[3mm]
Note that, in Case 1.2, the invertibility of all limit operators of $\op(\bA)$ as well as the uniform boundedness of the norms of their inverses follows already from the invertibility of $A$. \\[2mm]
Now consider {\bf Case 2}, i.e., suppose that none of the $h(n)$ belongs to $\cup v_k Y_k^{-1}$. For $n \in \sN$, let $r_n$ stand for the smallest non-negative integer such that there is a $k_n \in \sN$ with $h(n) \in v_{k_n} (\partial_\Omega Y_{k_n})^{-1} \Omega_{r_n}$. Consequently,
\[
h(n) \not\in v_{k_n} (\partial_\Omega Y_{k_n})^{-1} \Omega_{r_n - 1} \quad \mbox{for all} \; n.
\]
Again we set $r^* := \liminf r_n$ and distinguish two cases. \\[2mm]
{\bf Case 2.1: $r^*$ is finite.} We proceed as in Case 1.1 and find a subsequence of $h$ (denoted by $h$ again) and an element $w_* \in \Gamma$ such that $h(n) \in v_{k_n} (\partial_\Omega Y_{k_n})^{-1} w_*$. Since the inclusion $h(n) \in v_{k_n} Y_{k_n}^{-1}$ in (\ref{e180209.1}) had not been used in Case 1.1 we can continue exactly as in that case to obtain that Theorem \ref{t180209.10} and its corollary hold verbatim in the case at hand, too. \\[2mm]
{\bf Case 2.2: $r^*$ is infinite.} As in Case 1.2, we choose the sequence $(r_n)$ as strongly monotonically increasing. Then we have
\begin{equation} \label{e190209.9}
h(n) \not\in v_k Y_k^{-1} \quad \mbox{for all} \; k, \, n,
\end{equation}
\begin{equation} \label{e190209.10}
h(n) \not\in v_k (\partial_\Omega Y_k)^{-1} \Omega_{r_n -1} \quad \mbox{for all} \; k, \, n.
\end{equation}
We claim that these two facts imply that
\begin{equation} \label{e190209.12}
\Omega_{r_n -1} \cap Y_k v_k^{-1} h(n) = \emptyset \quad \mbox{for all} \; k, \, n.
\end{equation}
Indeed, from (\ref{e190209.9}) we conclude that $e \not\in Y_k v_k^{-1} h(n)$. Thus, for each $k$ and $n$, $\Omega_{r_n -1}$ contains points from the complement of $Y_k v_k^{-1} h(n)$, for instance the point $e$. Suppose that $\Omega_{r_n -1}$ also contains points in $Y_k v_k^{-1} h(n)$. Then the arguments from Case 1.2 imply that $\Omega_{r_n -1}$ contains points in the $\Omega$-boundary of $Y_k v_k^{-1} h(n)$. But (\ref{e190209.10}) implies that
$\Omega_{r_n -1} \cap (\partial_\Omega Y_k) v_k^{-1} h(n) = \emptyset$. Thus, $\Omega_{r_n -1}$ is completely located in the  complement of $Y_k v_k^{-1} h(n)$, whence (\ref{e190209.12}).

Since the operator $R_{h(n)}^{-1} R_{v_k} A_k R_{v_k}^{-1} R_{h(n)}$ lives on $\im P_{Y_k v_k^{-1} h(n)}$, we obtain from (\ref{e190209.12})
\begin{eqnarray*}
R_{h(n)}^{-1} (\op(\bA) + P_{\Gamma^\prime}) R_{h(n)}
& = & \sum_{k \ge 1} R_{h(n)}^{-1} R_{v_k} A_k R_{v_k}^{-1} R_{h(n)} (I - P_{\Omega_{r_n - 1}}) \\
&& \quad + \; R_{h(n)}^{-1} P_{\Gamma^\prime} R_{h(n)} (I - P_{\Omega_{r_n - 1}}) + P_{\Omega_{r_n - 1}}.
\end{eqnarray*}
The first two summands on the right-hand side of this equality tend strongly to zero as $n \to \infty$, whereas the third one tends strongly to the identity. Thus, the identity operator is the only limit operator of $\op(\bA) + P_{\Gamma^\prime}$ in Case 2.2. The following theorem summarizes the results from Cases 1.1 - 2.2.
\begin{theo} \label{t190209.13}
Let $\bA \in \cS_\cY ({\sf BDO}(\Gamma))$ and $A := \mbox{\rm s-lim} \, A_n P_{Y_n}$. Then the limit operators of $\op(\bA) + P_{\Gamma^\prime}$ are the identity operator $I$, all shifts $R_{v^*}^{-1} A R_{v^*}$ of the operator $A$, all limit operators of $A$, and all operators of the form
\[
\mbox{\rm s-lim} \, R_{w_*}^{-1} R_{\eta_{k_n}}^{-1} A_{k_n} R_{\eta_{k_n}} R_{w_*} + (P_{\Gamma^\prime})_g
\]
with a suitable subsequence $g$ of $h$ and with elements $\eta_{k_n} \in (\partial_\Omega Y_{k_n})^{-1}$ and $w_* \in \Gamma$.
\end{theo}
Combining this theorem with Theorems \ref{t160209.6} $(b)$ , \ref{t1.1} and \ref{t120309.1} we arrive at the following stability results.
\begin{theo} \label{t190209.14}
Let $\Gamma$ be a finitely generated exact discrete group, and let $(A_n) \in \cS_\cY ({\sf BDO}(\Gamma))$. The sequence $(A_n)$ is stable if and only if the operator 
\[
A := \mbox{\rm s-lim} \, A_n P_{Y_n}
\]
and all operators of the form
\[
\mbox{\rm s-lim} \,  R_{\eta_{k_n}}^{-1} A_{k_n} R_{\eta_{k_n}} + R_{w_*} (P_{\Gamma^\prime})_g R_{w_*}^{-1}
\]
with a suitable subsequence $g$ of $h$ and with elements $\eta_{k_n} \in (\partial_\Omega Y_{k_n})^{-1}$ and $w_* \in \Gamma$ are invertible and if the norms of their inverses are uniformly bounded.
\end{theo}
\begin{theo} \label{c190209.15}
Let $\Gamma$ be an exact discrete group, and let $A \in {\sf BDO}(\Gamma)$. The sequence $\bA = (P_{Y_n} A P_{Y_n})$ is stable if and only if the operator $A$ and all operators
\[
P_{\cY^{(h)}} A_g P_{\cY^{(h)}} : \im P_{\cY^{(h)}} \to \im P_{\cY^{(h)}}
\]
where $h$ is a sequence such that the limit $(\ref{e190209.1})$ exists and $\cY^{(h)}$ is as in $(\ref{e190209.16})$ and where $g$ is in $\sigma_{op, \, h}(A)$ are invertible and if the norms of their inverses are uniformly bounded.
\end{theo}
\begin{theo} \label{t120309.1a}
Let $\Gamma$ be a finitely generated discrete and exact group with sub-exponential growth which possesses at least one non-cyclic element, and let $A$ be a band operator on $l^2(\Gamma)$. Then the sequence $\bA = (P_{Y_n} A P_{Y_n})$ is stable if and only if the operators mentioned in the previous theorem are invertible.
\end{theo}
There are special sequences $\cY = (Y_n)$ and $\eta : \sN \to \Gamma$ for which the existence of the limit (\ref{e190209.16}) can be guaranteed. Let again $\Omega_n$ refer to the set of all products of at most $n$ elements of $\Omega$ and set $\Omega_0 := \{e\}$. A sequence $(\nu_n)$ in $\Gamma$ is called a {\em geodesic path} (with respect to $\Omega$) if there is a sequence $(w_n)$ in $\Omega \setminus \{e\}$ such that $\nu_n = w_1 w_2 \ldots w_n$ and $\nu_n \in \Omega_n \setminus \Omega_{n-1}$ for each $n \ge 1$. Note that this condition implies that each $\nu_n$ is in the {\em right $\Omega$-boundary of} $\Omega_n$, which is the set of all $w \in \Omega_n$ for which $w \Omega$ is not a subset of $\Omega_n$.

We will see now that the $\lim \Omega_n \eta_n$ exists if $\eta$ is an inverse geodesic path, i.e., if $\eta_n = \nu_n^{-1}$ for a geodesic path $\nu$.
\begin{lemma} \label{l190209.17}
Let $(w_n)_{n \ge 1}$ be a sequence in $\Omega$ and set $\eta_n := w_n^{-1} w_{n-1}^{-1} \ldots w_1^{-1}$ for $n \ge 1$. Then the strong limit $\mbox{\rm s-lim} \, P_{\Omega_n \eta_n}$ exists, and
\begin{equation} \label{e190209.18}
\mbox{\rm s-lim} \, P_{\Omega_n \eta_n} = P_{\cup_{n \ge 1} \Omega_n \eta_n}.
\end{equation}
\end{lemma}
{\bf Proof.} For $n \ge 1$, one has $\Omega_n \eta_n = \Omega_n w_{n+1} w_{n+1}^{-1} w_n^{-1} \ldots w_1^{-1} \subseteq \Omega_{n+1} \eta_{n+1}$. These inclusions imply the existence of the strong limit and the equality (\ref{e190209.18}). \hfill \qed \\[3mm]
The natural question arises whether every sequence $\eta : \sN \to \Gamma$ for which the set limit (\ref{e190209.16}) exists has a subsequence which is a subsequence of an inverse geodesic path. If the answer is affirmative, then it would prove sufficient to consider strong limits with respect to inverse geodesic paths in Theorem \ref{t190209.14} and its corollary. Under some conditions, this question was answered in \cite{Roc12} for commutative groups $\Gamma$ and for the free (non-commutative) groups $\sF_N$ with $N$ generators.
{\small Authors' addresses: \\[2mm]
Vladimir S. Rabinovich, Instituto Politechnico Nacional, ESIME
Zacatenco, Ed. 1, 2-do piso, Av. IPN, Mexico, D.F., 07738, Mexico.\\
Email: rabinov@maya.esimez.ipn.mx \\[1mm]
Steffen Roch, Technische Universit\"at Darmstadt, Fachbereich
Mathematik, Schlossgartenstrasse 7, D 64289 Darmstadt, FRG.\\
E-mail: roch@mathematik.tu-darmstadt.de}

\begin{thebibliography}{11}
%
%
\bibitem{Ada1}
{\sc T. Adachi}, A note on the F{\o}lner condition for amenability. -- Nagoya Math. J. {\bf 131}(1993), 67 -- 74.
\bibitem{Bla1}
{\sc B. Blackadar}, $K$-Theory for Operator Algebras. -- M. S. R. I. Monographs, Vol. 5, Springer-Verlag, Berlin, New York 1986.
\bibitem{BroO}
{\sc N. P. Brown, N. Ozawa}, $C^*$-Algebras and Finite-Dimensional Approximations. -- Graduate Studies Math. {\bf 88}, Amer. Math. Soc., Providence, R. I., 2008.
\bibitem{Dav1}
{\sc K. R. Davidson}, $C^*$-Algebras by Example. -- Fields Institute Monographs Vol. 6, Providence, R. I., 1996.
\bibitem{HRS2}
{\sc R. Hagen, S. Roch, B. Silbermann}, $C^*$-Algebras and Numerical Analysis. -- Marcel Dekker, Inc., New York, Basel 2001.
\bibitem{Lin1}
{\sc M. Lindner}, The finite section method and stable subsequences. -- Preprint 15/2008 TU Chemnitz, to appear in J. Appl. Num. Math.
\bibitem{Ped1}
{\sc G. K. Pedersen}, $C^*$-Algebras and Their Automorphism Groups. -- Academic Press, New York 1979.
\bibitem{RaR5}
{\sc V. S. Rabinovich, S. Roch}, Fredholm properties of band-dominated operators on periodic discrete structures. -- Complex Anal. Oper. Theory
{\bf 2}(2008), 4, 637 -- 681.
\bibitem{RRS1}
{\sc V. S. Rabinovich, S. Roch, B. Silbermann}, Fredholm theory and finite section method for band-dominated operators. -- Integral Equations Oper. Theory {\bf 30}(1998), 452 -- 495.
\bibitem{RRS2}
{\sc V. S. Rabinovich, S. Roch, B. Silbermann}, Band-dominated operators with operator-valued coefficients, their Fredholm properties and finite sections. -- Integral Eq. Oper. Theory {\bf 40}(2001), 3, 342 -- 381.
\bibitem{RRSB}
{\sc V. S. Rabinovich, S. Roch, B. Silbermann}, Limit Operators and their Applications in Operator Theory. -- Oper. Theory: Adv. Appl. {\bf 150}, Birkh\"{a}user, Basel 2004.
\bibitem{RRS6}
{\sc V. S. Rabinovich, S. Roch, B. Silbermann}, On finite sections of band-dominated operators. -- Operator Algebras, Operator Theory and Applications (Eds. M. A. Bastos, I. Gohberg, A. B. Lebre, F.-O. Speck), Oper. Theory: Adv. Appl. {\bf 181}, Birkh\"{a}user, Basel 2008, 385 -- 391.
\bibitem{Roc10}
{\sc S. Roch}, Finite sections of band-dominated operators. -- Memoirs AMS Vol. {\bf 191}, 895, Providence, R.I., 2008.
\bibitem{Roc11}
{\sc S. Roch}, Spatial discretization of Cuntz algebras. -- Houston Math. J., to appear.
\bibitem{Roc12}
{\sc S. Roch}, Spatial discretization of restricted group algebras. -- Preprint TU Darmstadt 2596, December 2009, 29 pages, submitted to Operators and Matrices.
\bibitem{Roe1}
{\sc J. Roe}, Lectures on Coarse Geometry. -- Univ. Lecture Ser. {\bf 31}, Amer. Math. Soc., Providence, R. I., 2003.
\bibitem{Roe2}
{\sc J. Roe}, Band-dominated Fredholm operators on discrete groups. --  Integral Equations Oper. Theory {\bf 51}(2005), 3, 411 -- 416.
%
\end{thebibliography}
\end{document}